\documentclass[12pt]{article}
\newcommand{\qed}{\nolinebreak \hfill{$\Box $} \par\vspace{0.5\parskip}\setcounter{equation}{0}}

\usepackage{txfonts}

\newcommand{\rea}{\hspace{2pt}\hbox{R\hspace{-13pt} I\hspace{7pt}}}

\newcommand{\natu}{\hspace{2pt}\hbox{N\hspace{-14pt} I\hspace{7pt}}}

\usepackage{url}
\usepackage{txfonts}
\usepackage{graphicx}
\usepackage{pdfpages}

\begin{document}

\begin{center}

{\large \bf An Algebraic Proof of Weierstrass's Approximation Theorem}
\vskip .5cm

{\it Jos\'e M. Gonz\'alez Barrios, Alberto Contreras-Crist\'an and Patricia I. Romero-Mares}

Instituto de Investigaciones en Matem\'aticas Aplicadas y en Sistemas.

Universidad Nacional Aut\'onoma de M\'exico
\end{center}

\begin{center}
    {\large\bf Abstract}
\end{center}

{\it In this paper we use the Vandermonde matrices and their
properties to give a new proof of the classical result of Karl
Weierstrass about the approximation of continuous functions $f$
on closed intervals, using a sequence of polynomials.
The proof solves  linear systems of equations using that
the Vandermonde matrices have always non zero determinants,
when the entries of the power series of the rows of
the matrix are all different.

We provide several examples, and we also use our method
to observe that the sequence of polynomials that we construct
algebraically approaches the Taylor series of a function $f$
which is infinitely differentiable.}

\noindent {\bf 1.-  Introduction} Karl Weierstrass has been
considered  one of the best
mathematicians in history. Many mathematicians tried to define the
concept of continuity including Augustin Louis Cauchy,
who provided a definition of continuity of functions.
However, the clearest and most accepted epsilon-delta
definition was given by Weierstrass, who also gave the
definition of uniformly continuous function, which is the
 basis of very important results in Mathematics.

Even more, Weierstrass gave the first example of a continuous function
which is not differentiable at any point, against the universal belief
that a continuous function had to be differentiable 
``almost everywhere". It is important to observe, that many mathematicians
did not like this result, for example Henri Poincar\'e called this
function an ``evil function'' which does not have any utility besides
being a counterexample, even when Weierstrass's result provided a very important
point in the creation of other areas in Mathematics such as Fractal Theory,
see \cite{Se}.

Weierstrass at the age of 70, also proved, using the concept of
uniformly continuous functions, that a continuous function 
defined on a closed interval
can be approximated by a sequence of polynomials. This result
had a great impact on the mathematical community.
In fact, Weierstrass published his first Theorem at age 39, and his
most important results appeared after in his later years. This
proof was received as a major result in mathematics,  and several
mathematicians proposed alternative proofs or extended substantially
the result, see for example,  \cite{B} who used what is called
the Bernstein basis polynomials, or \cite{S} who generalized 
extensively the result of
Weierstrass. This paper is an hommage to Karl Weierstrass.

There are several videos and pages that include proofs of the approximation theorem of
Weierstrass and giving references of this result. Among them 
\cite{emis}, \cite{wikiw} and \cite{people}. Most proofs use a family of
polynomials, integrals and convolutions to approximate $f$ a continuous
function defined on a closed interval. In this paper we will
use a simpler method based on the Vandermonde matrices to give
a nice approximation of $f$ even when $f$ is not a differentiable
function. The result is  applied in several examples that
prove that this technique is useful even for a not so  large degree
$n$ of the approximating polynomials.

\vskip 1cm
\noindent{\bf 2.- Main results}

Suppose that we have a continuous function $f:[a,b]\rightarrow \rea$,
 where $-\infty <a<b<\infty$. Then using Weierstrass's definition, we know 
 that $f$ is uniformly continuous. In fact, since $[a,b]$ is a compact interval in $\rea$, we can assume that $[a,b]=[0,1]$, by moving
 and rescaling the function 
from $[a,b]$ to $[0,1]$, by taking 
 the bijection $h(x)=\frac{x-a}{b-a}$ and defining
 $g:[0,1]\rightarrow \rea$ by $g(y)=f(a+y(b-a))$.
 The main properties of $f$ are preserved by this well-known
 transformation.
In the following Theorem we will use this result.

\noindent {\bf Theorem 2.1} {\it  Let $f:[0,1]\rightarrow \rea$ be a continuous function. Then
 there exists $N\in \natu$  and a polynomial $P_N(x)=a_N x^N+\cdots+ a_1 x^1+ a_0$, such that defining $g(t)=f(t)-P_N(t)$ 
for every $t\in [0,1]$, we have that
for every $\epsilon>0$  there exists $\delta>0$ satisfying that if $x,y \in [0,1]$ and
$|x-y|< \delta$ then 

\begin{equation}  |g(x)-g(y)| < \epsilon, \end{equation}
that is the polynomial $P_N$ approximates uniformly the function $f$.}   

\noindent {\bf Proof:} Since $f:[0,1]\rightarrow \rea$ is continuous and $[0,1]$ is a compact  set, 
$f$ is uniformly continuous, that is, for every $\epsilon>0$ there exists $\delta_1>0$, such that for every
$x,y\in [0,1]$ with $|x-y|<\delta_1$ we have that 

\begin{equation}\label{forf} |f(x)-f(y)|< \frac{\epsilon}{2}. \end{equation}

\noindent If $f$ is a polynomial of degree $n$, then $f$ is continuous on
$[0,1]$ and also uniformly continuous. So, we can define
$P_{n+j}(x)=f(x)$ for every $j\in\natu\cup\{0\}$, and there is nothing to prove.

\noindent For general continuous $f:[0,1]\rightarrow \rea$, 
let $N_0$ be the smallest integer such that $\frac{1}{N_0}< \delta_1$. Let
$$ x_0=0 <x_1=\frac{1}{N_0}<x_2=\frac{2}{N_0}<\cdots <x_{N_0}=\frac{N_0}{N_0}=1. $$
Find the polynomial $P_{N_0}$ defined by $\{x_0, x_1,\ldots ,x_{N_0} \}\subset [0,1]$ of degree $N_0$ such that
$P_{N_0}(j/N_0)=f(j/N_0)$ for every $j\in \{0,1,\ldots, N_0\}$, the existence
of this polynomial is proved below. Then

\begin{equation}\label{pn0}
P_{N_0}(x)= a_{N_0} x^{N_0}+ a_{N_0-1} x^{N_0-1}+\cdots +a_1 x^1 +a_0.
\end{equation}
Since $P_{N_0}$ is continuous on $[0,1]$ then it also uniformly continuous. 

Let $k\in \natu\cup\{0\}$  the set of nonnegative integers, and define

\begin{equation}\label{bk}
    N_k= 2^k \cdot N_0 \,\,\mbox{and}\,\,
    B_k=\{ 0=0/N_k, 1/N_k, 2/N_k,\ldots, (N_k-1)/N_k, N_k/N_k=1\},
\end{equation}
clearly $\{B_k\}_{k=0}^{\infty}$ is a family of partitions of the
closed interval $[0,1]$, such that 
$B_0\subset B_1\subset B_2\subset\cdots$ and if we define
$B=\cup_{k=0}^{\infty} B_k$, then $B$ is a countable dense subset of
$[0,1]$. For each $k\in\natu\cup\{0\}$ let

$$ x_0=0 <x_1=\frac{1}{N_k}<x_2=\frac{2}{N_k}<\cdots <
x_{N_k-1}=\frac{N_{k-1}}{N_k}<x_{N_k}=\frac{N_K}{N_K}=1. $$ 
As in (\ref{pn0}), and using (\ref{bk}) we define  polynomials $P_{N_k}$,
such that for every $j/N_k\in B_k$, $f(j/N_k)=P_{N_k}(j/N_k)$
where

\begin{equation}\label{pnk}
P_{N_k}(x)= a_{N_k} x^{N_k}+ a_{N_k-1} x^{N_k-1}+\cdots +a_1 x^1 +a_0.
\end{equation}
The existence of these polynomials $P_{N_k}$ in (\ref{pnk}) is proved below. 

It is clear that $B_k$ is a uniform partition of $[0,1]$, which
 approaches  a countable dense set of $[0,1]$. Hence, there must exist
 a large $k\in\natu\cup\{0\}$ such that the sup distance between
 $f$ and $P_{N_k}$ can be smaller than $\epsilon/2$ for the $\epsilon$
 given in (\ref{forf}). Observe that for every $x\in[0,1]$, we have
  $\lim_{k\rightarrow\infty} P_{N_k}(x)=f(x)$ using continuity. \qed

Recall that if $n\in\natu$ is an integer greater than $0$ and we consider $[a,b]$ a closed interval such that $-\infty <a <b <\infty$,
then $[a,b]$ is a nonempty compact set, then we know that a {\bf bounded polynomial of degree} $n$  is a function
$f:[a,b]\rightarrow \rea$ such that $f(x)=a_n x^n + a_{n-1} x^{n-1} + \cdots + a_1 x^1 +a_0$ for every $x\in [a,b]$,
where  $\{a_0, a_1, \ldots , a_{n-1}, a_n \}\subset \rea$ is the set of the $n+1$ coefficients of the polynomial and $a_n\not =0$.

\setcounter{equation}{5}
Now, let $m\in\natu$ another integer and define a subset of $m+1$ different points in the compact set $[a,b]$ such that
\begin{equation} \label{m+1points}
a \leq x_1 < x_2 < \cdots < x_m <x_{m+1}\leq b.
\end{equation}
Define two square matrices of order $(m+1)\times (m+1)$  given by

\begin{equation} \label{Amatrix}
 A = \left( 
\begin{array}{ccccc}
x_1^m & x_1^{m-1} & \cdots & x_1^1 & 1 \nonumber \\
x_2^m & x_2^{m-1} & \cdots & x_2^1 & 1 \nonumber \\
\vdots& \vdots    & \vdots & \vdots& \vdots \nonumber \\
x_m^m & x_m^{m-1} & \cdots & x_m^1 & 1 \nonumber \\
x_{m+1}^m & x_{m+1}^{m-1} & \cdots & x_{m+1}^1 & 1 
\end{array} \right)\,\,\mbox{and}\,\,
 B = \left( 
\begin{array}{ccccc}
1 & x_1^1 & \cdots & x_1^{m-1} & x_1^{m} \nonumber \\
1 & x_2^1 & \cdots & x_2^{m-1} & x_2^{m} \nonumber \\
\vdots& \vdots    & \vdots & \vdots& \vdots \nonumber \\
1 & x_m^1 & \cdots & x_m^{m-1} & x_m^{m} \nonumber \\
1 & x_{m+1}^1& \cdots & x_{m+1}^{m-1} & x_{m+1}^{m} 
\end{array} \right).
\end{equation}
The right-hand side matrix $B$ is known as the Vandermonde
matrix of order $m+1\times m+1$, named after A. T. Vandermonde,
a french mathematician of the seventeenth century, it is also
known that $B$ is always an invertible matrix, since its
determinant $\mbox{Det}(B)=\Pi_{1\leq i<j\leq m+1}(x_i-x_j)$
is always different from $0$, when all the $x_i's$
are different.  We also know that the exchange of two columns
of the matrix $B$ multiplies the determinant by $-1$, see \cite{H97} or
\cite{wikid}.
If we exchange the first and last columns of $B$, then we exchange
the second and the penultimate columns of $B$, and so on, then
$\mbox{Det(B)} = \pm \mbox{Det}(A)$ depending on the number
of columns of $B$, if this number  is even then we have the plus sign,
otherwise we have the minus sign. Hence, the determinant
of the matrix $B$ is not zero if and only if the determinant of $A$
in (\ref{Amatrix}) is not zero. 

\noindent The proof of $Det(B)\not =0$ is given in Appendix 1.

\vskip 1cm

\noindent Now we study by cases the relation between $n+1$ 
the number of coefficients $\{ a_n,a_{n-1},\ldots ,a_{1},a_{0}\}$
in a polynomial of degree $n$,  and $m+1$ the size of the matrix $A$ defined
in (\ref{Amatrix}).

\noindent In the case that $f$ is a polynomial of degree $n=m$, it is important to observe that  for every
set of $n+1$ points given by
\begin{equation}\label{coeff}
a\leq x_1 <x_2<\cdots x_n<x_{n+1} \leq b,\end{equation}
the solution of the linear equation $A \cdot \underline{a} = \underline{b}$,
where 
$$\underline{a}=(a_n, a_{n-1}, \ldots ,a_1, a_0)^t$$
is the transpose vector of the coefficients of $f$, and
$$\underline{b} = (b_1=f(x_1), b_2=f(x_2), \ldots , b_n=f(x_n),
b_{n+1}=f(x_{n+1}))^t$$
is the transpose vector of the evaluations of $f$, then the solution
of the equation, which exists by Remark A in Appendix 1 is given by
$$\underline{a} = A^{-1}\cdot \underline{b},$$
independently of the different points $x_1, x_2, \ldots ,x_n, x_{n+1}$
satisfying (\ref{coeff})

\noindent Second case, assume that $m>n$ and $f$ is a polynomial of degree
$n$. Since $m>n$, there exists an integer $k\geq 2$, such that 
$m+1=n+k$, so let
\begin{equation}
a\leq x_1<x_2<\cdots <x_n< \cdots <x_{n+k}\leq b,
\end{equation}
In this case we have that
$$\underline{a}=(a_{n+k},\ldots,a_{n+1}, a_n,\ldots, a_1, a_0)^t$$
and 
$$\underline{b}=(b_1, b_2, \ldots , b_{n+1}, \ldots ,b_{n+k})^t.$$
Here again the matrix $A$ of size $(n+k)\times(n+k)$ has determinant
different from zero, so $A$ is invertible. But since $f$ is a
polynomial of degree $n$, then the coefficients
$a_{n+1}, a_{n+2}, \ldots , a_{n+k}$ need to be zero, and  in fact
the solution of $A\cdot \underline{a}=\underline{b}$ given by
$\underline{a} = A^{-1}\cdot \underline{b}$, features such result.

\noindent Third case, assume that $m<n$, where $n$
is the degree of the polynomial, in this case we will
find a polynomial of degree $n-j$ for some $1\leq j<n$.
It is  very important to observe that if we have
$n-j+1$ points given by
$$a\leq x_1<x_2<\cdots <x_{n-j}<x_{n-j+1}\leq b $$
then by Remark A in Appendix 1, the linear system defined  by these
$n-j+1$ points and the matrix $A$ of order $(n-j+1)\times(n-j+1)$
given in (\ref{Amatrix}) has a non zero determinant and the
solution $\underline{a} = A^{-1}\cdot  \underline{b}$ exists and is well
defined, but the solution depends strongly in the selection
of the points 
$a\leq x_1<x_2<\cdots <x_{n-j}<x_{n-j+1}\leq b$. This
statement will be clarified in the examples below.

It is also important to observe that the solution of the linear system
$\underline{a} = A^{-1} \cdot \underline{b}$,
where $A$ is defined as in equation (\ref{Amatrix}), has the same solution
if we permute the rows of the matrix $A$ and we use the same
permutation in the vector $\underline{b}$ getting $\underline{b_p}$.
It is known that if we have an invertible matrix $A$ 
of order $m\times m$, and we rearrange 
the rows using a permutation $p$ of $\{ 1,2,\ldots ,m\}$,
getting a matrix $A_p$,
then its  inverse $A_p^{-1}$ is simply the matrix $A^{-1}$
after permuting    columns using $p$, and
$A_p^{-1}\cdot \underline{b_p}=\underline{a}$.
We give a very simple example assume that
$f(x)=-x^2+3x+5=a_1 x^2+a_2 x+a_3$ and $f$ is defined on $[-1,1]$.
Let $-1=x_1<x_2=0<x_3=1$, then $f(x_1)=1$, $f(x_2)=5$ and $f(x_3)=7$
define $A$ as in (\ref{Amatrix}), then
$$A=\left(\begin{array}{ccc}
1 & -1 & 1\nonumber\\
0 &  0 & 1 \nonumber\\
1 & 1 & 1 
\end{array}\right),$$
whose determinant is $\mbox{Det}(A)=-2$,
and its inverse is
$$A^{-1}=\left(\begin{array}{ccc}
1/2 & -1 & 1/2\nonumber\\
-1/2 &  0 & 1/2 \nonumber\\
0 & 1 & 0 
\end{array}\right).$$
Now define $\underline{b}=(1, 5 , 7)^{t}$, and the solution
of the linear system $\underline{a}=A^{-1}\cdot \underline{b}$ is

\begin{equation}\left(\begin{array}{c}
a_1\\
a_2\\
a_3
\end{array}\right)=
\left(\begin{array}{ccc}
1/2 & -1 & 1/2\\
-1/2 &  0 & 1/2\\
0 & 1 & 0 
\end{array}\right)\cdot
\left(\begin{array}{c}
1\\
5\\
7
\end{array}\right)=
\left(\begin{array}{c}
-1 \nonumber\\
3 \nonumber\\
5
\end{array}\right)
\end{equation}
Let $p:\{1,2,3\}\rightarrow \{1,2,3\}$ the permutation
such $p(1)=3, p(2)=2$ and $p(3)=1$, then the matrix $A_p$ associated
to $p$ is given by
$$A_p=\left(\begin{array}{ccc}
1 & 1 & 1\nonumber\\
0 &  0 & 1 \nonumber\\
1 & -1 & 1 
\end{array}\right),$$
whose determinant is $\mbox{Det}(A)=2$,
and its inverse is
$$A_{p}^{-1}=\left(\begin{array}{ccc}
1/2 & -1 & 1/2\nonumber\\
1/2 &  0 & -1/2 \nonumber\\
0 & 1 & 0 
\end{array}\right).$$
Here we observe that the columns of $A^{-1}$ have been permuted
using the permutation $p$, that is, the first column
of $A^{-1}$ has moved to the third column of $A_p^{-1}$,
and the third column of $A^{-1}$ has moved to the first column of
$A_p^{_1}$ and the second in both inverses has not moved.
Now we permute the entries of $\underline{b}$, to obtain
$\underline{b_p}=(7,5,1)^t$. In this case
we have that
\begin{equation}
A_p^{-1}\cdot \underline{b_p} = \left(\begin{array}{ccc}
1/2 & -1 & 1/2\\
1/2 &  0 & -1/2 \\
0 & 1 & 0 
\end{array}\right)\cdot
\left(\begin{array}{c}
7\\
5\\
1
\end{array}\right)=
\left(\begin{array}{c}
-1 \nonumber\\
3 \nonumber\\
5,
\end{array}\right)
\end{equation}
and $(-1,3,5)=(a_1,a_2,a_3)=\underline{a}$.

Now we will provide some examples when the continuous function
is a polynomial of some finite degree $n$.

\noindent {\bf Example 2.3} Let $f(x)=x+3$ defined on $[-1,4]$, then $f$ is a polynomial of degree $n=1$ defined on a compact interval.

First, let us assume that $-1 \leq x_1=-1<x_2=1<4$, then
$f(x_1)=f(-1)=2$ and $f(x_2)=f(1)=4$.
Let $\underline{a} =(a_1 ,a_0)^t$ and $\underline{b}=(f(x_1), f(x_2))^t$.
Then $A$ the matrix associated
to $x_1$ and $x_2$ and the equation $A\cdot\underline{a}=\underline{b}$
is given by
$$ A = \left( \begin{array}{cc}
-1 & 1 \nonumber\\
1 & 1
\end{array}\right),$$
whose inverse is given by

$$ A^{-1} = \left( \begin{array}{cc}
-1/2 & 1/2 \nonumber\\
1/2 & 1/2
\end{array}\right).$$
So, the solution is given by

\begin{equation}\label{gr1=11}
\left(\begin{array}{c}
a_1 \\
a_0
\end{array}\right) =\left( \begin{array}{cc}
-1/2 & 1/2 \nonumber\\
1/2 & 1/2
\end{array}\right) \cdot 
\left(\begin{array}{c}
2 \\
4
\end{array}\right)
=\left(\begin{array}{c}
1 \\ 
3
\end{array}\right).
\end{equation}

Second,  let us assume that $-1 < x_1=0<x_2=2<4$, then
$f(x_1)=f(0)=3$ and $f(x_2)=f(2)=5$.
Let $\underline{a} =(a_1 ,a_0)^t$ and $\underline{b}=(f(x_1), f(x_2))^t$.
Then $A$ the matrix associated
to $x_1$ and $x_2$ and the equation $A\cdot\underline{a}=\underline{b}$
is given by
$$ A = \left( \begin{array}{cc}
0 & 1 \nonumber\\
2 & 1
\end{array}\right),$$
whose inverse is given by

$$ A^{-1} = \left( \begin{array}{cc}
-1/2 & 1/2 \nonumber\\
1    & 0
\end{array}\right).$$
So, the solution is given by

\begin{equation}\label{gr1=12}
\left(\begin{array}{c}
a_1 \\
a_0
\end{array}\right) =\left( \begin{array}{cc}
-1/2 & 1/2 \nonumber\\
 1   & 0
\end{array}\right) \cdot 
\left(\begin{array}{c}
3 \\
5
\end{array}\right)
=\left(\begin{array}{c}
1 \\ 
3
\end{array}\right).
\end{equation}
As observed above, both solutions in (\ref{gr1=11}) and (\ref{gr1=12})
are equal, because the degree of the polynomial is $n=1$.

Third, let us assume now that we have three points
$-1\leq x_1=-1<x_2=0<x_3=1$, then $f(x_1)=2$, $f(x_2)=3$ and
$f(x_3)=4$ in this case we have that
$\underline{a}=(a_2, a_1, a_0)^t$ and
$\underline{b}=(2,3,4)^t$ and the matrix $A$ is given by

$$ A = \left( \begin{array}{ccc}
1 & -1 & 1 \nonumber\\
0 & 0  & 1  \nonumber\\
1 & 1  & 1
\end{array}\right),$$
whose inverse is given by

$$ A^{-1} = \left( \begin{array}{ccc}
1/2  & -1 & 1/2 \nonumber\\
-1/2 & 0  & 1/2 \nonumber\\
0    &  1 & 0
\end{array}\right).$$
So, the solution is given by

\begin{equation}\label{gr1=+}
\left(\begin{array}{c}
a_2 \\
a_1 \\
a_0
\end{array}\right) =
\left( \begin{array}{ccc}
1/2 & -1 &  1/2 \nonumber\\
-1/2& 0  & 1/2 \nonumber\\
0   & 1  & 0
\end{array}\right) \cdot 
\left(\begin{array}{c}
2 \\
3 \\
4
\end{array}\right)
=\left(\begin{array}{c}
0 \\
1 \\ 
3
\end{array}\right).
\end{equation}
Therefore, the coefficient of $x_2$ is $0$, which is the correct number.

Finally, let us assume that we have four points
$-1\leq x_1=-1<x_2=0<x_3=1<x_4=2<4$, then
$f(x_1)=2, f(x_2)=3, f(x_3)=4$ and $f(x_4)=5$, so we have that
$\underline{a}=(a_3, a_2,a_1, a_0)^t$
and $\underline{b}=(2,3,4,5)^t$ and the matrix $A$ is given by

$$ A = \left( \begin{array}{cccc}
-1 & 1 & -1 & 1 \nonumber\\
0  & 0 & 0  & 1 \nonumber\\
1  & 1 & 1  & 1 \nonumber \\
8  & 4 & 2  & 1
\end{array}\right),$$
whose inverse is given by

$$ A^{-1} = \left( \begin{array}{cccc}
-1/6 & 1/2  & -1/2 & 1/6 \nonumber\\
1/2  & -1   & 1/2  & 0   \nonumber\\
-1/3 & -1/2 &  1   &-1/6 \nonumber\\
0    &  1   & 0    & 0
\end{array}\right).$$
So, the solution is given by

\begin{equation}\label{gr1=++}
\left(\begin{array}{c}
a_3 \\
a_2 \\
a_1 \\
a_0 
\end{array}\right) =
\left( \begin{array}{cccc}
-1/6 & 3/6 & -3/6 & 1/6 \nonumber\\
1/2  &-2/2 & 1/2  & 0   \nonumber\\
-2/8 & -3/6& 6/6  & -1/6\nonumber \\
0    & 1   & 0    & 0
\end{array}\right) \cdot 
\left(\begin{array}{c}
2 \\
3 \\
4 \\
5
\end{array}\right)
=\left(\begin{array}{c}
0 \\
0 \\
1 \\ 
3
\end{array}\right).
\end{equation}
Here, as expected the coefficients $a_3=a_2=0$.

\vskip .5cm
\noindent {\bf Example 2.4} Let us consider the polynomial defined
on $[-1,1]$ such that
$f(x)=x^4 +2$. First take two points $-1\leq x_1=-1 <x_2=1\leq 1$,
 then $f(x_1)=3=f(x_2)$. In this case we have that
 $\underline{a}=(a_1,a_0)^t$ and $\underline{b}=(3,3)^t$
 the matrix $A$
 is given by

 $$ A = \left( \begin{array}{cc}
-1 & 1 \nonumber\\
1 & 1
\end{array}\right),$$
whose inverse is given by

$$ A^{-1} = \left( \begin{array}{cc}
-1/2 & 1/2 \nonumber\\
1/2 & 1/2
\end{array}\right).$$
So, the solution is given by

\begin{equation}\label{gr1=}
\left(\begin{array}{c}
a_1 \\
a_0
\end{array}\right) =\left( \begin{array}{cc}
-1/2 & 1/2 \nonumber\\
1/2 & 1/2
\end{array}\right) \cdot 
\left(\begin{array}{c}
3 \\
3
\end{array}\right)
=\left(\begin{array}{c}
0 \\ 
3
\end{array}\right).
\end{equation}
In this case a polynomial of degree $1$ is in fact a polynomial of degree
$0$, that is $P_1(x)=0x+3$.

Second, now we consider three points
$-1\leq x_1=-1<x_2=0<x_3=1\leq 1$, then
$f(-1)=3, f(0)=2$ and $f(1)=3$. Here,
$\underline{a}=(a_2,a_1,a_0)^t$ and $\underline{b}=(3,2,3)^t$,
and the corresponding matrix $A$ is given by

$$ A = \left( \begin{array}{ccc}
1 & -1 & 1 \nonumber\\
0 & 0  & 1  \nonumber\\
1 & 1  & 1
\end{array}\right),$$
whose inverse is given by

$$ A^{-1} = \left( \begin{array}{ccc}
1/2  & -1 & 1/2 \nonumber\\
-1/2 & 0  & 1/2 \nonumber\\
0    &  1 & 0
\end{array}\right).$$
So, the solution is given by

\begin{equation}\label{gr1=+mas}
\left(\begin{array}{c}
a_2 \\
a_1 \\
a_0
\end{array}\right) =
\left( \begin{array}{ccc}
1/2 & -1 &  1/2 \nonumber\\
-1/2& 0  & 1/2 \nonumber\\
0   & 1  & 0
\end{array}\right) \cdot 
\left(\begin{array}{c}
3 \\
2 \\
3
\end{array}\right)
=\left(\begin{array}{c}
1 \\
0 \\ 
2
\end{array}\right).
\end{equation}
So, a polynomial of order two approximating $f$
 is $P_2(x)=x^2+2$.

 Third, selecting four points  
 $-1\leq x_1=-1<x_2=-1/3<x_3=1/3<x_4=1\leq 1$, then
 $f(-1)=f(1)=3$ and $f(-1/3)=f(1/3)=163/81$, then
 $\underline{a}=(a_3,a_2,a_1,a_0)^t$ and
 $\underline{b}= (3,163/81,163/81,3)^t$ and the matrix $A$ is

 $$ A = \left( \begin{array}{cccc}
-1    & 1  & -1  & 1 \nonumber\\
-1/27 &1/9 & -1/3& 1 \nonumber\\
1/27  &1/9 & 1/3 & 1 \nonumber \\
1     & 1  & 1   & 1
\end{array}\right),$$
whose inverse is given by

$$ A^{-1} = \left( \begin{array}{cccc}
-9/16 & 27/16  & -27/16 & 9/16 \nonumber\\
9/16  & -9/16   & -9/16  & 9/16   \nonumber\\
1/16 & -27/16 &  27/16   &-1/16 \nonumber\\
-1/16 &  9/16  & 9/16 & -1/16
\end{array}\right).$$
So, the solution is given by

\begin{equation}\label{gr1=+-}
\left(\begin{array}{c}
a_3 \\
a_2 \\
a_1 \\
a_0
\end{array}\right) =
\left( \begin{array}{cccc}
-9/16 & 27/16  & -27/16 & 9/16 \nonumber\\
9/16  & -9/16   & -9/16  & 9/16   \nonumber\\
1/16 & -27/16 &  27/16   &-1/16 \nonumber\\
-1/16 &  9/16  & 9/16 & -1/16
\end{array}\right)\cdot
\left(\begin{array}{c}
3 \\
163/81 \\
163/81  \\
3
\end{array}\right)
=\left(\begin{array}{c}
0    \\
10/9 \\
0    \\ 
17/9
\end{array}\right).
\end{equation}
Hence, we get the polynomial $P_3(x)= 0x^3+10/9 x^2+0x+17/9=10/9 x^2+17/9$,
which satisfies the values of $f(-1)=3, f(-1/3)=163/81, f(1/3)=163/81$ and $f(1)=3$.

Fourth, let us assume that we take five points given by
$-1\leq x_1=-1<x_2=-1/2<x_3=0< x_4=1/2<x_5=1\leq 1$,
then $f(x_1)=3, f(x_2)=33/16, f(x_3)=2, f(x_4)=33/16$ and $f(x_5)=3$. Here,
$\underline{a} =(a_4,a_3,a_2,a_1,a_0)^t$ and
$\underline{b}=(3,33/16,2,33/16,3)^t$, and the matrix $A$ is defined as

$$ A = \left( \begin{array}{ccccc}
1   &-1  & 1   & -1  & 1   \nonumber\\
1/16&-1/8& 1/4 &-1/2 & 1    \nonumber \\
0   &0   & 0   & 0   & 1      \nonumber \\
1/16&1/8 &1/4  &1/2  & 1  \nonumber \\
1   &1   & 1   & 1   & 1  \nonumber
\end{array}\right)
$$
and its inverse is given by
$$ A^{-1} = \left( \begin{array}{ccccc}
2/3   & -8/3 &12/3  & -8/3 & 2/3  \nonumber\\
-2/3  & 4/3  & 0    &-4/3  & 2/3  \nonumber \\
-1/6  & 16/6 &-5    & 16/6 & -1/6 \nonumber \\
1/6   & -8/6 & 0    & 8/6  &-1/6  \nonumber \\
0     &0     & 1    &   0  &  0  \nonumber
\end{array}\right).
$$
The solution of the system is

\begin{equation} \left(\begin{array}{c}
a_4 \nonumber\\
a_3 \nonumber\\
a_2 \nonumber\\
a_1 \nonumber\\
a_0 \nonumber

\end{array}\right)= \left( \begin{array}{ccccc}
2/3   & -8/3 &12/3  & -8/3 & 2/3  \nonumber\\
-2/3  & 4/3  & 0    &-4/3  & 2/3  \nonumber \\
-1/6  & 16/6 &-5    & 16/6 & -1/6 \nonumber \\
1/6   & -8/6 & 0    & 8/6  &-1/6  \nonumber \\
0     &0     & 1    &   0  &  0  \nonumber
\end{array}\right) \cdot\left(\begin{array}{c}
3    \nonumber \\
33/16\nonumber \\
2    \nonumber \\
33/16\nonumber \\
3    \nonumber
\end{array}\right)=\left(\begin{array}{c}
1   \nonumber\\
0   \nonumber\\
0   \nonumber\\
0  \nonumber\\
2   \nonumber
\end{array}\right)
\end{equation}
Hence, the solution is the polynomial $P_4(x)=x^4+2$.  If we add more points
then we get that the coefficients of the terms with
exponents higher than or equal to $5$ are
zeros.

\noindent {\bf Example 2.5} 
Suppose that we have a polynomial but we do not know its degree and it is
defined for every  $x\in [0,4]$ assume that we are given
three points $0\leq x_1=1< x_2=2<x_3=3\leq 4$, and we also get their values under $f$, $f(1)=4$, $f(2)=26$ and
$f(3)=86$, in this case $\underline{a}=(a_2,a_1,a_0)^t$ and
$\underline{b}=(4,26,86)^t$ and $A$ the matrix associated is given by

$$ A = \left( \begin{array}{ccc}
1   & 1   & 1  \nonumber \\
4   & 2   & 1  \nonumber \\
9   & 3   & 1
\end{array}\right)
$$
and its inverse is given by
$$ A^{-1} = \left( \begin{array}{ccc}
1/2 & -1  & 1/2  \nonumber \\
-5/2& 4   & -3/2 \nonumber \\
3   & -3  & 1
\end{array}\right).
$$
The solution of the system is 

\begin{equation}\label{2sol}\left(\begin{array}{c}
a_2 \nonumber \\
a_1 \nonumber \\
a_0  \nonumber
\end{array}\right)= \left( \begin{array}{ccc}
1/2 & -1  & 1/2  \nonumber \\
-5/2& 4   & -3/2 \nonumber \\
3   & -3  & 1   \nonumber
\end{array}\right) \cdot \left(\begin{array}{c}
4   \nonumber\\
26  \nonumber\\
86  \nonumber
\end{array}\right)=\left(\begin{array}{c}
19  \nonumber \\
-35 \nonumber\\
20  
\end{array}\right)
\end{equation}
Hence the solution is $P_2(x) =19x^2-35x+20$, and we have that $P_2(1)=4$, $P_2(2)=26$ and $P_2(3)=86$.

Now assume that we are given other three points $0\leq x_1=2<x_2=3<x_3=4\leq 4$ and their values under $f$,
$f(2)=26$, $f(3)=86$ and $f(4)=202$, in this case $\underline{a}=(a_2,a_1,a_0)^t$ and
$\underline{b}=(26,86,202)^t$ and $A$ the matrix associated is given by

$$ A = \left( \begin{array}{ccc}
4   & 2   & 1  \nonumber \\
9   & 3   & 1  \nonumber \\
16  & 4   & 1
\end{array}\right)
$$
and its inverse is given by
$$ A^{-1} = \left( \begin{array}{ccc}
1/2 & -1  & 1/2  \nonumber \\
-7/2& 6   & -5/2 \nonumber \\
6   & -8  & 3
\end{array}\right).
$$
The solution of the system is

\begin{equation} \left(\begin{array}{c}
a_2 \nonumber\\
a_1 \nonumber\\
a_0  \nonumber
\end{array}\right)= \left( \begin{array}{ccc}
1/2 & -1  & 1/2  \nonumber \\
-7/2& 6   & -5/2 \nonumber \\
6   & -8  & 3  \nonumber
\end{array}\right) \cdot \left(\begin{array}{c}
26  \nonumber \\
86  \nonumber \\
202 \nonumber
\end{array}\right)=\left(\begin{array}{c}
28  \nonumber\\
-80 \nonumber\\
74 \nonumber
\end{array}\right)
\end{equation}
Here the solution is $P_2(x) =28x^2-80x+74$, and we have that $P_2(2)=26$, $P_2(3)=86$ and $P_2(4)=202$,
and this solution is not equal to the solution in (\ref{2sol}). Therefore, the degree of $f$ is not $n=2$.

Now we take the four points given above $ 0\leq x_1=1<x_2=2<x_3=3<x_4=4\leq 4$, then we know that $f(1)=4, f(2)=26, f(3)=86$ and $f(4)=202$.
Here, $\underline{a}=(a_3, a_2, a_1, a_0)^t$ and $\underline{b}=(4,26,86,202)^t$ and the matrix $A$ associated  with these vectors is given by

$$ A = \left( \begin{array}{cccc}
1   & 1   & 1   & 1    \nonumber \\
8   & 4   & 2   & 1      \nonumber \\
27  & 9   & 3   & 1  \nonumber \\
64  & 16  & 4   & 1
\end{array}\right)
$$
and its inverse is given by
$$ A^{-1} = \left( \begin{array}{cccc}
-1/6 & 1/2 & -1/2  & 1/6  \nonumber \\
3/2  & -4  & 7/2   & -1   \nonumber \\
-13/3& 19/2& -7    &11/6  \nonumber \\
4    & -6  &   4   & -1
\end{array}\right).
$$
The solution of the system is

\begin{equation} \left(\begin{array}{c}
a_3 \nonumber\\
a_2 \nonumber\\
a_1 \nonumber\\
a_0 \nonumber
\end{array}\right)= \left( \begin{array}{cccc}
-1/6 & 1/2 & -1/2  & 1/6  \nonumber \\
3/2  & -4  & 7/2   & -1   \nonumber \\
-13/3& 19/2& -7    &11/6  \nonumber \\
4    & -6  &   4   & -1   \nonumber
\end{array}\right) \cdot\left(\begin{array}{c}
4   \nonumber \\
26  \nonumber \\
86  \nonumber \\
202 \nonumber
\end{array}\right)=\left(\begin{array}{c}
3   \nonumber\\
1   \nonumber\\
-2  \nonumber\\
2   \nonumber
\end{array}\right)
\end{equation}
Hence the polynomial of degree $n=3$ is given by $P_3(x)=3x^3+x^2 -2x+2$, which satisfies
that $P_3(1)=4, P_3(2)=26, P_3(3)=86$ and $P_3(4)=202$. In fact, it is easy to see that
if we take the four points $0\leq x_1=0<x_2=1<x_3=2<x_4=3\leq 4$, we obtain again the same polynomial $P_3$.

Finally, we are told that $f(0)=2$, then we take five points given by
$0\leq x_1=0<x_2=1<x_3=2<x_4=3<x_5=4\leq 4$, let
$\underline{a}=(a_4,a_3,a_2,a_1,a_0)^t$ and $\underline{b}=(2,4,26,86,202)^t$ in this
case the matrix $A$ is given by

$$ A = \left( \begin{array}{ccccc}
0   &0   & 0   & 0   & 1   \nonumber\\
1   &1   & 1   & 1   & 1    \nonumber \\
16  &8   & 4   & 2   & 1      \nonumber \\
81  &27  & 9   & 3   & 1  \nonumber \\
256 &64  & 16  & 4   & 1  \nonumber
\end{array}\right)
$$
and its inverse is given by
$$ A^{-1} = \left( \begin{array}{ccccc}
1/24  & -1/6 & 1/4  & -1/6 & 1/24 \nonumber\\
-5/12 & 3/2  & -2   & 7/6  & -1/4 \nonumber \\
35/24 &-13/3 &19/4  & -7/3 & 11/24\nonumber \\
-25/12& 4    & -3   & 4/3  &-1/4  \nonumber \\
1     &0     & 0    &   0  &  0  \nonumber
\end{array}\right).
$$
The solution of the system is

\begin{equation} \left(\begin{array}{c}
a_4 \nonumber\\
a_3 \nonumber\\
a_2 \nonumber\\
a_1 \nonumber\\
a_0 \nonumber

\end{array}\right)= \left( \begin{array}{ccccc}
1/24  & -1/6 & 1/4  & -1/6 & 1/24 \nonumber\\
-5/12 & 3/2  & -2   & 7/6  & -1/4 \nonumber \\
35/24 &-13/3 &19/4  & -7/3 & 11/24\nonumber \\
-25/12& 4    & -3   & 4/3  &-1/4  \nonumber \\
1     &0     & 0    &   0  &  0  \nonumber
\end{array}\right) \cdot\left(\begin{array}{c}
2   \nonumber \\
4   \nonumber \\
26  \nonumber \\
86  \nonumber \\
202 \nonumber
\end{array}\right)=\left(\begin{array}{c}
0   \nonumber\\
3   \nonumber\\
1   \nonumber\\
-2  \nonumber\\
2   \nonumber
\end{array}\right)
\end{equation}
Therefore, the polynomial we obtain is $P_4(x) =0x^4+3x^3+x^2-2x+2=
P_3(x)$, where $P_3$ is given above. Hence, now we know that
$f$ is a polynomial of degree $n=3$.

\noindent {\bf Example 2.6}
Now we give a polynomial approximating the absolute value function
defined on $[-1,1]$, Here we use a polynomial of degree $n=18$
that approaches nicely the absolute value function, even when this function
is not differentiable at $x=0$. See Figure 1.

\newpage
\phantom{x}
\begin{figure}[h!]
\includepdf[pages=-, scale=.8]{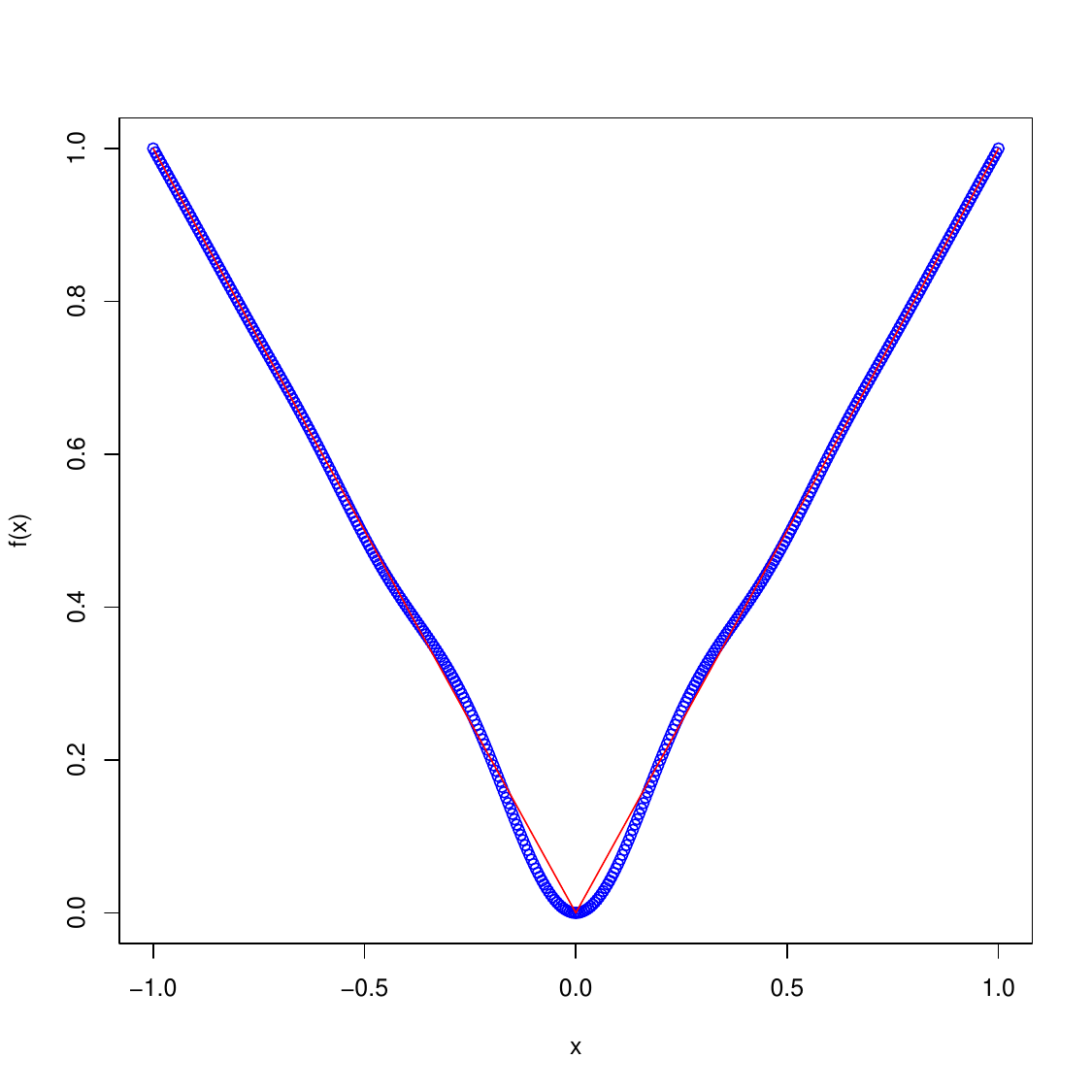}
\vskip 18cm
\caption{Approximation of the absolute value function}
\end{figure}
\vskip 8cm


The set of points selected on $[-1,1]$ were
$x_1=-1,x_2=-.96,x_3=-.92,x_4=-.88,x_5=-.84,x_6=-.8,x_7=-.6,
x_8=-.4,x_9=-.2,$
$x_{10}=0,x_{11}=.2,x_{12}=.4,x_{13}=.6,x_{14}=.8,x_{15}=.84$
$x_{16}=,88,x_{17}=.92,x_{18}=.98,x_{19}=1$
Using the solution of the equation with absolute values
we obtain the polynomial
\begin{eqnarray}
f(x)&=&158.266x^{18}-86,539X^{16}+2031.475x^{14}-2669.014x^{12}
\nonumber\\
& &+2143.315x^{10}-1077.715x^8+334.136x^6-61.044x^4+6.971x^2.
\end{eqnarray}

\noindent {\bf Example 2.7} Finally, we will find a polynomial that
approaches  the sine function in $[-\pi,\pi]$. In order to
facilitate the evaluation we will transform the original onto
$[0,1]$, that is, we will approximate the function

$$ g(y)= \sin(2\pi y - \pi) \quad \mbox{for every}\quad y\in [0,1].$$
We will use two polynomials of degrees $n=6$ and $n=8$.
For the case $n=6$ we will use the points
$x_1=0, x_2=1/6, x_3=2/6, x_4=3/6, x_5=4/6, x_6=5/6, x_7=1$.
Using the linear equation with the above points to generate the
matrix $A$, we obtained the following polynomial
$$  f(x)=5.57 e-11x^6+56.118 x^5-140.296 x^4+101.324 x^3-11.691x^2-5.456 x.$$

\begin{center}
\begin{figure}[h]
\includegraphics[width=15cm]{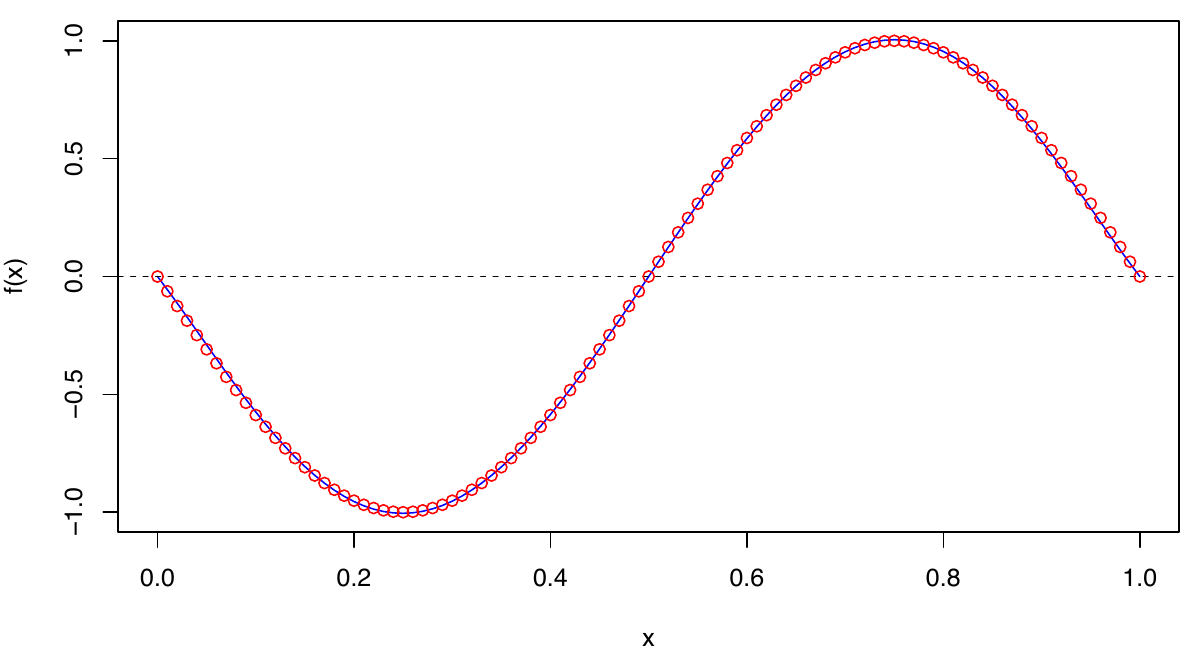}
\vskip -1cm
\caption{Graph of the polynomial for the sine function of degree n=6.}
\end{figure}
\end{center}

\
\begin{center}
\begin{figure}[h]
\includegraphics[width=15cm]{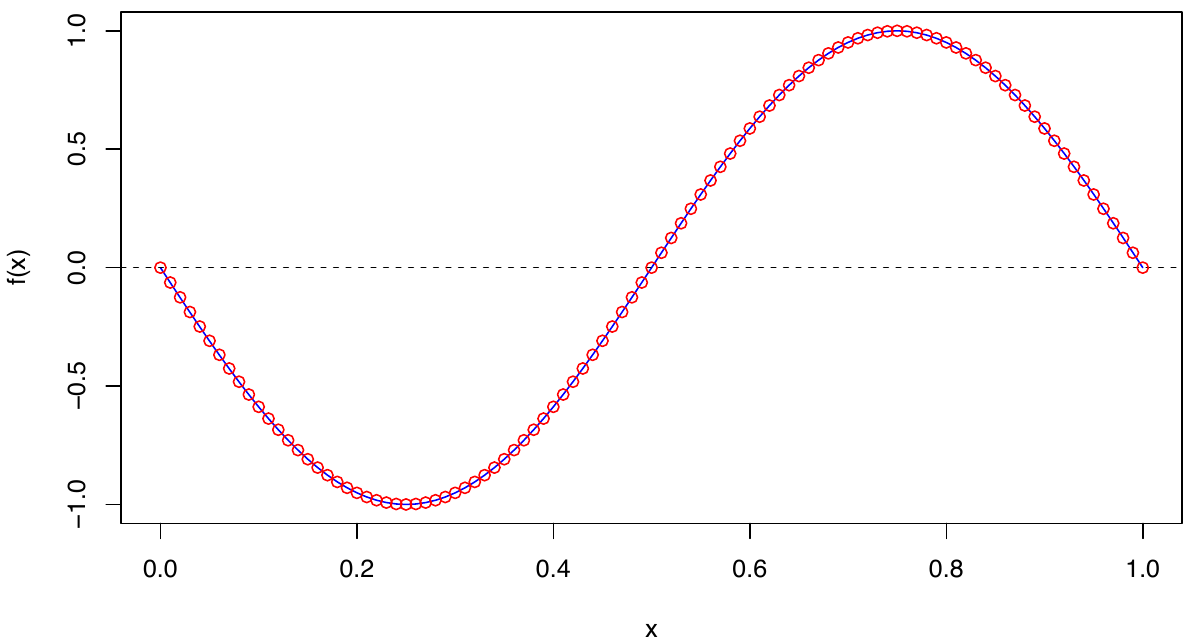}
\vskip -1cm
\caption{ Graph of the polynomial for the sine function of degree n=8.}
\end{figure}
\end{center}

The graph of this approximation  together with the 
graph of the sine function is given above.

In the case $n=8$ we used the points
$x_1=0, x_2=1/8, x_3=2/8, x_4=3/8, x_5=4/8, x_6=5/8, x_7=6/8, x_8=7/8$
and $x_9=1$. Using the same procedure as above, we find the 
Vandermonde matrix $A$ of the linear equation, and solving this
equation we obtain the polynomial

$$f(x)=-2.75e-12x^8 -59.143x^7+207x^6-231.305x^5+60.764x^4 +27.381x^3-
1.644x^2-6.361x.$$
The graph of this approximation together with the 
graph of the sine function is given above.

Now we will try to approximate the Taylor's series of the
function $\sin(x)$ around $x=0$ on the interval $[-\pi,\pi]$,
using the polynomials generated by the linear system of equations.
It is well known that the Taylor series of the sine function
around $x=0$ is given by
\begin{equation}\label{TaySin}
\sin(x)=\sum_{k=0}^{\infty}\frac{\sin^{k}(0)}{k!}x^k=
\sum_{k=0}^{\infty}\frac{(-1)^k \cdot x^{2k+1}}{(2k+1)!}.
\end{equation}
Since the sine function is an odd function, we would expect
that using our approximations the even powers of the sine
function would not appear in the approximating polynomial.
In fact, this is the case when we try the degrees of the 
approximating polynomials
of order $n=4, n=6, n=8$ and $n=10$ when we take the partitions of 
$[-\pi,\pi]$, $\Pi_4=\{-\pi,-\pi/2,0,\pi/2,\pi\}$,

\noindent $\Pi_6=\{-\pi,-2\pi/3, -\pi/3,0,\pi/3,2\pi/3,\pi\}$,
$\Pi_8=\{-\pi,-3\pi/4,\ldots,0,\ldots, 3\pi/4,\pi\}$ and

\noindent $\Pi_{10}=\{-\pi,-4\pi/5,\ldots,0,\ldots ,4\pi/5, \pi\}.$
Solving the linear systems, we obtain the following polynomials

\begin{equation}\label{p4p6}
P_4(x)=-\frac{-8}{3\pi^3} x^3+\frac{8}{3\pi} x  \quad\mbox{and}\quad
P_6(x)=\frac{81\sqrt{3}}{80\pi^5}x^5
-\frac{45}{16\pi^3} x^3+\frac{9\sqrt{3}}{5/\pi} x,
\end{equation}

\begin{equation}\label{p8}
P_8(x)=-\frac{2048(-7+5\sqrt{2})}{315\pi^7} x^7+
\frac{128(-26+19\sqrt{2})}{45\pi^5} x^5
-\frac{8(-169+140\sqrt{2})}{45\pi^2} x^3
+\frac{8(-21+44\sqrt{2})}{105\pi} x
\end{equation}

\begin{eqnarray}\label{p10}
P_{10}(x)&=&\frac{390625\left(34\sqrt{5-\sqrt{5}}-21\sqrt{5+\sqrt{5}}\right)}
{142152\sqrt{2}\,\pi^9}* x^9\nonumber\\
& &-\,\frac{15625\left(326\sqrt{5-\sqrt{5}}-201(\sqrt{5+\sqrt{5}}\right)}
{24192\sqrt{2}\,\pi^7}* x^7\nonumber\\
& &+\,\frac{625\left(1786\sqrt{5-\sqrt{5}}-1089(\sqrt{5+\sqrt{5}}\right)}
{6912\sqrt{2}\,\pi^5}* x^5\nonumber\\
& & -\,\frac{25\left(67576\sqrt{5-\sqrt{5}}-37821(\sqrt{5+\sqrt{5}}\right)}
{36288\sqrt{2}\,\pi^3}* x^3\nonumber\\
& &+\, \frac{25\left(83\sqrt{5-\sqrt{5}}-18(\sqrt{5+\sqrt{5}}\right)}
{504\sqrt{2}\,\pi}* x
\end{eqnarray}
In  Table 1 we compare the values obtained by the coefficients
of the polynomials
(\ref{p4p6}), (\ref{p8}) and (\ref{p10}) and the coefficients of the Taylor
series of the sine function given in (\ref{TaySin}).
As we can note the approximation given by the polynomials of degrees
$4,6,8$ and $10$ provide a very accurate value to the coefficients
in the Taylor expansion.

\begin{center}
\mbox{Table 1.- Approximation of the Coefficients of the Taylor Series of the Sine function}
\begin{tabular}{|l|c|c|c|c|c|}
\hline
\mbox{Value of} $k$ & 1 & 3 & 5 & 7 & 9 \\
\mbox{value of} $1/k!$ & 1 & 0.166666 &0.008333 & 0.0001984 &0.00000275\\
\mbox{approx. given by} $P_4$ & 0.848826 &0.086040 & - & - & -  \\
\mbox{approx. given by} $P_6$ & 0.992392 &0.157109 &0.005730 & - & - \\
\mbox{approx. given by} $P_8$ & 0.999805 &0.166216&0.008087&0.0001529 & - \\
\mbox{approx. given by} $P_{10}$&0.999997& 0.166655&0.008323&0.0001948&0.00000225\\
\hline
\end{tabular}
\end{center}
This example also shows that our approximation polynomials provide
a tool to approximate the values of the derivatives of the function
$f$ at $x_0=0$ in the Taylor series without differentiating.

\noindent {\bf Example 2.8}
Our last example is also related to Taylor series, in the case that
the radius of convergence of the series is bounded. We know that
 the function $f(x)=\ln(x+1)$ is well defined for $x>-1$, but
 $\lim_{x\rightarrow -1}f(x)=-\infty$ and we also know that the
 Taylor series, also known as the Maclaurin series of $f(x)$ around
 $ x=0$ is given by
 \begin{equation}\label{log}
 \ln(1+x)=\sum_{n=1}^{\infty} \frac{(-1)^{n-1}}{n} x^n
 \quad\mbox{for every}\quad |x|<1.
 \end{equation}
 Of course the convergence of this series has problems as
 $x$ approaches from the right to $-1$.

 We propose to make approximations using polynomials
 of the Taylor series of $f$ for 
 $x\in[-3/4,3/4]$ taking uniform partitions of this interval. We give only the approximations
 for degrees, $n=4$, $n=6$, $n=8$ and $n=10$. As above, we only 
 give in Table 2 the approximations for $k\leq 5$.

 As we can note the approximation given by the polynomials of degrees
$4,6,8$ and $10$ provide a very accurate value to the coefficients
in the Taylor expansion.

\begin{center}
\mbox{Table 2.- Approximation of the Coefficients of the Taylor Series of the  function $\ln(1+x)$}
\begin{tabular}{|l|c|c|c|c|c|}
\hline
\mbox{Value of} $k$ & 1 & 2 & 3 & 4 & 5 \\
\mbox{value of} $(-1)^{k-1}/k$ & 1 & -0.5 &0.333333 & -0.25 &0.02\\
\mbox{approx. given by} $P_4$ & 0.969277 &-0.453517 & 0.583103 & -0.464548 & -  \\
\mbox{approx. given by} $P_6$ & 1.003036 &-0.502726 &0.269673 & -0.192895 &0.450513 \\
\mbox{approx. given by} $P_8$ & 0.999703 &-0.499727&0.345072&-0.260784 &0.096422 \\
\mbox{approx. given by} $P_{10}$&1.000028&-0.500026&0.331491&-0.248281&0.228492\\
\hline
\end{tabular}
\end{center}
For the values of $k$ between $6$ and $10$ the approximations are not so
good because we are approaching the value of $x$  to $-1$.

\noindent 
\begin{center}{\bf Appendix 1}\end{center}

In the following Remark, we will use two  known properties of the determinant
of a square matrix, they are:

\noindent {\bf 1.-} Adding or subtracting a multiple of one column to another column
does not change the determinant of a matrix

\noindent{\bf 2.-} Multiplying a row of a matrix by a number $c$ multiplies the
determinant by the number $c$.

\noindent {\bf Remark A} {\it Let $-\infty <a <b <\infty$, let $n$ be an integer greater than or equal to $1$ and take
$a \leq x_1 < x_2 < \cdots < x_n \leq b$ as
in (\ref{m+1points}) with $m+1=n$ and define the square matrix $A$ as in equation (\ref{Amatrix}). Then the determinant 
of the matrix $A$ is not zero, so, there exists the inverse matrix of $A$, denoted by $A^{-1}$.}

Let $x_1, x_2, \ldots, x_{n}$ be $n$ different numbers, let $B$ a Vandermonde matrix 
of size $n\times n$ of the form

\begin{equation}\label{Vand}
B=V(x_1, x_2, \ldots, x_n)= \left(\begin{array}{cccccc}
1 & x_1 & x_1^2 &\cdots & x_1^{n-2} & x_1^{n-1} \nonumber\\
1 & x_2 & x_2^2 &\cdots & x_2^{n-2} & x_2^{n-1} \nonumber\\
\vdots &\vdots &\vdots &\vdots &\vdots &\vdots  \nonumber\\
1 & x_{n-1} & x_{n-1}^2 &\cdots & x_{n-1}^{n-2} & x_{n-1}^{n-1} \nonumber\\
1 & x_n & x_n^2 &\cdots & x_n^{n-2} & x_n^{n-1} 
\end{array}\right)
\end{equation}
We want to see by induction that 

\begin{equation}\label{detVan}
\mbox{Det}(B)= \Pi_{1\leq i< j\leq n} (x_j-x_i).
\end{equation}

Let $n=2$, then 

$$\mbox{Det}(V(x_1,x_2))=\mbox{Det}\left(\begin{array}{cc}
1 & x_1 \nonumber\\
1 & x_2 \nonumber\\
\end{array}\right) = x_2-x_1,
$$
which is (\ref{detVan}) for $n=2$.

\noindent Let us assume the Induction Hypothesis, that is,  
$\mbox{Det}(V(x_1,x_2,\ldots,x_{n-1}))=\Pi_{1\leq i<j\leq n-1}(x_j-x_i)\not= 0$,
for $n-1$ different numbers $x_1,x_2, \ldots, , x_{n-1}$.

\noindent Let $x_1,x_2, \ldots, , x_{n-1}, x_n $, $n$ different numbers and consider 
the matrix $B=V(x_1, x_2, \ldots, x_n)$ defined
on (\ref{Vand}). Let $C_1,C_2, \ldots ,C_{n-1}, C_n$ be the $n$ columns of $B$, and redefined the columns
$C_n =C_n -x_n\cdot C_{n-1}$, $C_{n-1}=C_{n-1}-x_n\cdot C_{n-2}$, 
$\cdots$, $C_3=C_3-x_n\cdot C_2$ and $C_2=C_2-x_n\cdot C_1$, using 1.-,
the determinant does not change and we get

\begin{equation}\label{detcol}
\mbox{Det}(B)= \mbox{Det}\left(\begin{array}{cccc}
1 & x_1-x_n & \cdots & x_1^{n-1}-x_1^{n-2}x_n \nonumber\\
1 & x_2-x_n & \cdots & x_2^{n-1}-x_2^{n-2}x_n \nonumber\\
\vdots & \vdots &  \vdots & \vdots \nonumber\\
1 & x_{n-1}-x_n  &\cdots & x_{n-1}^{n-1}-x_{n-1}^{n-2}x_n \nonumber\\
1 & x_n-x_n             &\cdots & x_n^{n-1}-x_n^{n-1}\nonumber\\
\end{array}\right) = 
\mbox{Det}\left(\begin{array}{cccc}
1 & x_1-x_n &\cdots & x_1^{n-1}-x_1^{n-2}x_n \nonumber\\
1 & x_2-x_n &\cdots & x_2^{n-1}-x_2^{n-2}x_n \nonumber\\
\vdots & \vdots & \vdots & \vdots \nonumber\\
1 & x_{n-1}-x_n &\cdots & x_{n-1}^{n-1}-x_{n-1}^{n-2}x_n \nonumber\\
1 & 0                               &\cdots & 0                              \nonumber\\
\end{array}\right)
\end{equation}
Expanding the determinant using the last row and the entry with row $n$ and column $1$ in (\ref{detcol}) 

\begin{equation}\label{detcol2}
\mbox{Det}(B)=(-1)^{n+1}\mbox{Det} \left(\begin{array}{cccc}
(x_1-x_n) & (x_1-x_n)\cdot x_1 & \cdots & (x_1-x_n)\cdot x_1^{n-2} \nonumber\\
(x_2-x_n) & (x_2-x_n)\cdot x_2 & \cdots & (x_2-x_n)\cdot x_2^{n-2} \nonumber\\
\vdots    & \vdots           & \vdots & \vdots  \nonumber\\
(x_{n-1}-x_n) & (x_{n-1}-x_n)\cdot x_{n-1} & \cdots & (x_{n-1}-x_n)\cdot x_{n-1}^{n-2} \nonumber\\
\end{array}\right)
\end{equation}
Factorizing from each row $(x_i-x_n)$ in (\ref{detcol2}),
using 2.- we obtain that the determinant is multiplied by all these
constants

\begin{equation}\label{detcol3}
\mbox{Det}(V(x_1,\ldots,x_n))=(-1)^{n+1}\Pi_{1\leq i\leq n-1} (x_i-x_n)\mbox{Det} \left(\begin{array}{ccccc}
1 & x_1 & x_1^2 & \cdots & x_1^{n-2} \nonumber\\
1 & x_2 & x_2^2 & \cdots & x_2^{n-2} \nonumber\\
\vdots & \vdots & \vdots &\vdots \vdots \nonumber\\
1 & x_{n-1} & x_{n-1}^2 & \cdots & x_{n-1}^{n-2} \nonumber\\
\end{array}\right)
\end{equation}
Now we use the induction hypothesis in (\ref{detcol3}) to finish the proof. \qed

\noindent {\bf Final Comments} 
In this paper, we observed that we can approximate a continuous
function $f$ defined on closed intervals $[a,b]$, or simply $[0,1]$, 
using  polynomials which are obtained solving a linear
system with invertible matrices, using a grid that is a partition
of the closed interval, in many instances a large uniform partition
is quite adequate when the function $f$ is not differentiable
in many points. The matrices $A$ and $A^{-1}$ only depend
on the points selected and not on the function $f$.

It is also remarkable to observe that if $f$ is infinitely 
differentiable we can approximate closely the Taylor series
of the function $f$ around a point $x_0$, obtaining
very close approximations of the values of the
$\{f^{(k)}(x_0)\}_{k=0}^{L}$, without differentiating,
when $L$ is not so small, by using polynomials with degrees a little
larger than $L$.

We are already working in the multivariate version of these results.

\noindent {\bf Acknowledgment} We want to thank our
colleague Sandra Palau
for valuable comments that greatly improve this paper.


\begin{thebibliography}{99}

\bibitem{B}
Bernstein, S.N. (1912/13). D\'emostration du Th\'eor\`em de
Weierstrass fond\'ee sur le calcul de probabilit\'es.
{\it Communications of the Kharkov Mathematical Society}. 
{\bf XIII}, p. 1-2.

\bibitem{CH}
Cheney, E. W. (2000).
{\it Introduction to Approximation Theory.} (2nd. Ed.).
Providence RI: AMS Chelsea Publ.
ISBN 978-0-8218-1374-4.


\bibitem{de}
de la Cerda, S. (2023). Polynomial Approximations to 
Continuous Functions. {\it The American Mathematical Monthly}.
{\bf 130}, (7), 655. doi:101080/00029890.2023.2206324.

\bibitem{H97}
Harville, D.A. (1997). {\it Matrix Algebra From a Statistician's Perspective.}
Springer-Verlag, New York. 


\bibitem{P}
Pinkus, A. (2013). Weierstrass and Approximation Theory.
{\it Journal of Approximation Theory}. {\bf 107}, (1):8
ISSN 0021-9045, OCLC 4638498762.


\bibitem{Se}
Sepulcre, J.M. (2016). {\it La gestaci\'on del An\'alisis moderno. Weierstrass.}
Colecci\'on Genios de la Matem\'aticas. 2017. Coleccionables, S.A.U.
RBA. Printed in Spain.

\bibitem{S}
Stone, M.H. (1948). The Generalized Weierstrass Approximation
Theorem. {\it Mathematics Magazine.} {\bf 21}, (4), 167--184.
doi:102307/30229750.

\bibitem{emis}
emis.de/journals/SAT/papers/99/4.pdf

\bibitem{wikid}
en.wikipedia.org/wiki/Determinant.


\bibitem{wikiw}
en.wikipedia.org/wiki/Stone-Weierstrass-theorem


\bibitem{people}
people.math.sc.edu/schep/weierstrass.pdf
\end{thebibliography}
\end{document}